\documentclass[journal]{IEEEtran}
\usepackage{cite}
\usepackage{amsmath,amsfonts,amssymb,amsthm}
\theoremstyle{definition}

\theoremstyle{plain}

\newtheorem{lemma}{Lemma}
\newtheorem{proposition}{Proposition}

\usepackage{graphicx}
\usepackage{mathrsfs}
\usepackage{booktabs}
\usepackage{multirow}
\usepackage{array}
\usepackage[ruled]{algorithm}
\usepackage{hyperref}
\hypersetup{colorlinks,%
	citecolor={red}, % Changer pour "black" avec natbib
	urlcolor={blue},
	linkcolor={blue},
	breaklinks={true},
	pagebackref={true},
	hyperindex={true},
}
\usepackage{url}
\usepackage{xcolor}
\usepackage{pgf,tikz}
\usetikzlibrary{arrows,shapes,positioning}
\usepackage{circuitikz}
\usepackage{subfig}
\usepackage{setspace}
\usepackage{algorithmic}
\floatname{algorithm}{Model}

\usepackage{enumitem}
\newlist{abbrv}{itemize}{1}
\setlist[abbrv,1]{label=,labelwidth=0.9in,align=parleft,noitemsep,leftmargin=!}
\newcommand{\R}{\mathbb{R}}
\newcommand{\C}{\mathbb{C}}

\newcommand{\Herm}{\mathbb{H}}
\newcommand{\G}{\mathscr{G}}

\newcommand{\E}{\mathcal{E}}
\newcommand{\K}{\mathcal{K}}
\newcommand{\cplx}[1]{\mathrm{#1}}
\newcommand{\rv}[1]{\boldsymbol{#1}}
\newcommand{\cv}[1]{\boldsymbol{\mathrm{#1}}}
\newcommand{\ub}[1]{\overline{#1}}
\newcommand{\lb}[1]{\underline{#1}}
\newcommand{\Node}{\mathcal{N}}
\newcommand{\Gen}{\mathcal{G}}
\newcommand{\Branch}{\mathcal{L}}
\newcommand{\Reseau}{\mathscr{P}}

\newcommand{\jc}{\cplx{j}}
\newcommand{\pl}{\lb{p}}
\newcommand{\ql}{\lb{q}}
\newcommand{\pu}{\ub{p}}
\newcommand{\qu}{\ub{q}}
\newcommand{\vl}{\lb{v}}
\newcommand{\vu}{\ub{v}}
\DeclareMathOperator{\subj}{s.t.}
\DeclareMathOperator{\rank}{rank}

\DeclareMathOperator{\re}{Re}
\DeclareMathOperator{\im}{Im}

\begin{document}
%
% paper title
% Titles are generally capitalized except for words such as a, an, and, as,
% at, but, by, for, in, nor, of, on, or, the, to and up, which are usually
% not capitalized unless they are the first or last word of the title.
% Linebreaks \\ can be used within to get better formatting as desired.
% Do not put math or special symbols in the title.
\title{Tight-and-Cheap Conic Relaxation for the \\ AC Optimal Power Flow Problem}
%
%
% author names and IEEE memberships
% note positions of commas and nonbreaking spaces ( ~ ) LaTeX will not break
% a structure at a ~ so this keeps an author's name from being broken across
% two lines.
% use \thanks{} to gain access to the first footnote area
% a separate \thanks must be used for each paragraph as LaTeX2e's \thanks
% was not built to handle multiple paragraphs
%

\author{Christian~Bingane,~\IEEEmembership{Student Member,~IEEE,}
        Miguel~F.~Anjos,~\IEEEmembership{Senior Member,~IEEE}
        and~S\'ebastien~Le~Digabel %,~\IEEEmembership{Life~Fellow,~IEEE}% <-this % stops a space
\thanks{The authors are with the Department of Mathematics and Industrial Engineering, Polytechnique Montreal, Montreal, Quebec, Canada H3C~3A7; and the GERAD research center, Montreal, Quebec, Canada H3T~2A7. E-mails: christian.bingane@polymtl.ca, anjos@stanfordalumni.org, sebastien.le-digabel@polymtl.ca.}% <-this % stops a space
\thanks{This research was supported by the NSERC-Hydro-Quebec-Schneider Electric Industrial Research Chair.}% <-this % stops a space
%\thanks{Manuscript received Mmmm dd, yyyy; revised Mmmm dd, yy.}
}

\maketitle

% As a general rule, do not put math, special symbols or citations
% in the abstract or keywords.
\begin{abstract}
The classical alternating current optimal power flow problem is highly nonconvex and generally hard to solve. Convex relaxations, in particular semidefinite, second-order cone, convex quadratic, and linear relaxations, have recently attracted significant interest. The semidefinite relaxation is the strongest among them and is exact for many cases. However, the computational efficiency for solving large-scale semidefinite optimization is lower than for second-order cone optimization. We propose a conic relaxation obtained by combining semidefinite optimization with the reformulation-linearization technique, commonly known as RLT. The proposed relaxation is stronger than the second-order cone relaxation and nearly as tight as the standard semidefinite relaxation. Computational experiments using standard test cases with up to 6515 buses show that the time to solve the new conic relaxation is up to one order of magnitude lower than for the chordal relaxation, a semidefinite relaxation technique that exploits the sparsity of power networks.
\end{abstract}

% Note that keywords are not normally used for peerreview papers.

% For peer review papers, you can put extra information on the cover
% page as needed:
% \ifCLASSOPTIONpeerreview
% \begin{center} \bfseries EDICS Category: 3-BBND \end{center}
% \fi
%
% For peerreview papers, this IEEEtran command inserts a page break and
% creates the second title. It will be ignored for other modes.
\IEEEpeerreviewmaketitle

\section*{Nomenclature}
\subsection{Notations}
\begin{abbrv}
\item[$\R$/$\C$] Set of real/complex numbers,
%\item[$\Sym^n$] Set of $n\times n$ real symmetric matrices
\item[$\Herm^n$] Set of $n\times n$ Hermitian matrices,
%\item[$\Herm^n_+$] Set of $n\times n$ positive semidefinite matrices,
\item[$\jc$] Imaginary unit,
\item[$a$/$\cplx{a}$] Real/complex number,
\item[$\rv{a}$/$\cv{a}$] Real/complex vector,
\item[$A$/$\cplx{A}$]  Real/complex matrix.
\end{abbrv}
\subsection{Operators}
\begin{abbrv}
\item[$\re (\cdot)$/$\im (\cdot)$]  Real/imaginary part operator,
\item[$(\cdot)^*$] Conjugate operator,
\item[$\left|\cdot\right|$]  Magnitude or cardinality set operator,
\item[$\angle (\cdot)$]  Phase operator,
\item[$(\cdot)^T$]  Transpose operator,
\item[$(\cdot)^H$]  Conjugate transpose operator,
\item[$\rank (\cdot)$]  Rank operator.
\end{abbrv}
\subsection{Input data}
\begin{abbrv}
\item[$\Reseau = (\Node, \Branch)$] Power network,
\item[$\Node$] Set of buses,
\item[$\Gen = \bigcup_{k\in\Node} \Gen_k$] Set of generators,
\item[$\Gen_k$] Set of generators connected to bus~$k$,
\item[$\Branch$] Set of branches,
\item[$p_{Dk}$/$q_{Dk}$] Active/reactive power demand at bus~$k$,
\item[$g_k'$/$b_k'$] Conductance/susceptance of shunt element at bus~$k$,
\item[$c_{g2}, c_{g1}, c_{g0}$] Generation cost coefficients of generator~$g$,
\item[$\cplx{y}_\ell^{-1} = r_{\ell} + \jc x_{\ell}$] Series impedance of branch~$\ell$,
\item[$b'_\ell$] Total shunt susceptance of branch~$\ell$,
\item[$\cplx{t}_{\ell}$] Turns ratio of branch~$\ell$.
\end{abbrv}
\subsection{Variables}
\begin{abbrv}
\item[$p_{Gg}$/$q_{Gg}$] Active/reactive power generation by generator~$g$,
\item[$\cplx{v}_{k}$] Complex (phasor) voltage at bus~$k$,
\item[$p_{f\ell}$/$q_{f\ell}$] Active/reactive power flow injected along branch~$\ell$ by its \emph{from} end,
\item[$p_{t\ell}$/$q_{t\ell}$] Active/reactive power flow injected along branch~$\ell$ by its \emph{to} end.
\end{abbrv}

\section{Introduction}
\IEEEPARstart{T}{he} optimal power flow (OPF) problem, introduced by Carpentier in 1962~\cite{ref2}, seeks to find a network operating point that optimizes an objective function such as generation cost subject to power flow equations and other operational constraints. A broad survey of the history of the problem and the related optimization methods appears in~\cite{cain, frank1, frank2}.

According to~\cite{frank1}, the general OPF problem may be modelled using linear, mixed-integer linear, nonlinear, or mixed integer nonlinear optimization. 
We focus on the nonlinear version, also called alternating current optimal power flow (ACOPF) problem. The ACOPF problem is nonconvex and NP-hard~\cite{verma2010,lehmann2016}. One way to tackle it is to use convex relaxations of the nonconvex constraints.

A conic optimization problem is a class of convex optimization problem that consists in optimizing a linear function over the intersection of an affine subspace and a convex cone. When the cone is the nonnegative orthant, the second-order cone, or the positive semidefinite matrices set, the conic optimization problem is a linear, a second-order cone or a semidefinite program respectively. A large theory can be found in~\cite{boydlivre} on convex optimization or in~\cite{sdplivre} on semidefinite optimization.

Since the ACOPF problem can be cast as a quadratically constrained quadratic program (QCQP), two principal conic relaxations have been proposed in the last decade: the second-order cone programming (SOCP) relaxation~\cite{jabr1} and the semidefinite programming (SDP) relaxation~\cite{bai}. These two relaxations offer several advantages. First, they can lead to global optimality. Second, because they are relaxations, they provide a bound on the global optimal value of the ACOPF problem. Third, if one of them is infeasible, then the ACOPF problem is infeasible.

We should note that the SDP relaxation is stronger than the SOCP relaxation but requires heavier computation. Therefore, a chordal relaxation was proposed in~\cite{jabr4} in order to exploit the fact that power networks are not densely connected, thus reducing data storage and increasing computation speed. A full literature review on these three relaxations can be found in~\cite{low1, low2}. Other convex relaxations have been developed in~\cite{hijazi2014, hijazi2016, josz2015, coffrin2016, mhanna2016}.

For radial networks, the SOCP relaxation is tantamount to the SDP relaxation. In this case, one would normally solve the first one rather than the second one due to the difference in computation time. For general meshed networks, it would be interesting to develop a relaxation as fast as the SOCP relaxation and as strong as the SDP relaxation. For example, three strong SOCP relaxations were developed in~\cite{kocuk2016} that are very close in quality to the SDP relaxation and are faster to solve. 

In this paper, we present a new conic relaxation that offers a favourable trade-off between the SOCP and the SDP relaxations for large-scale instances of ACOPF in terms of optimality gap and computation time. This relaxation is obtained through a combination of semidefinite optimization and the reformulation-linearization technique, known as RLT.

The remainder of this paper is organized as follows. In Section~\ref{sec1:formulation}, we define the mathematical model of the ACOPF problem (without loss of generality). In Section~\ref{sec1:conic}, we describe principal conic relaxations of the ACOPF problem, especially semidefinite and second-order cone relaxations. In Section~\ref{sec1:new}, we present the new conic relaxation, and we present computational results in Section~\ref{sec1:results}. Section~\ref{sec1:conclusion} concludes the paper.

\section{ACOPF: Formulation}\label{sec1:formulation}
Consider a typical power network $\Reseau = (\Node, \Branch)$, where $\Node = \{1,\ldots,n\}$ and $\Branch \subseteq \Node \times \Node$ denote respectively the set of buses and the set of branches (transmission lines, transformers and phase shifters). Each branch~$\ell \in \Branch$ has a \emph{from} end  $k$ (on the \emph{tap side}) and a \emph{to} end $m$ as modeled in~\cite{matpower}. We note $\ell = (k,m)$.  The ACOPF problem is given as:
\begin{subequations}\label{eq1:acopf}
\begin{equation}\label{eq1:objfun}
\min \sum_{g\in\Gen} c_{g2} p_{Gg}^2 + c_{g1} p_{Gg} + c_{g0}
\end{equation}
over variables $\rv{p}_{G}, \rv{q}_{G} \in \R^{|\Gen|}$, $\rv{p}_{f},\rv{q}_{f},\rv{p}_{t},\rv{q}_{t} \in \R^{|\Branch|}$, and $\cv{v}\in~\C^{|\Node|}$, subject to
\begin{itemize}
\item Power balance equations:
\begin{align}
&\sum_{g\in \Gen_k} p_{Gg} - p_{Dk} - g_k'\left|\cplx{v}_k\right|^2 = \nonumber\\
&\sum_{\ell=(k,m)\in\Branch} p_{f\ell} + \sum_{\ell=(m,k)\in\Branch} p_{t\ell}\;\forall k\in\Node,\label{eq1:kclp}\\
&\sum_{g\in \Gen_k} q_{Gg} - q_{Dk} + b_k'\left|\cplx{v}_k\right|^2 = \nonumber\\ &\sum_{\ell=(k,m)\in\Branch} q_{f\ell} + \sum_{\ell=(m,k)\in\Branch} q_{t\ell}\;\forall k\in\Node,\label{eq1:kclq}
\end{align}
\item Line flow equations:
\begin{align}
\frac{\cplx{v}_k}{\cplx{t}_{\ell}} &\left[\left(\jc \frac{b'_{\ell}}{2} + \cplx{y}_\ell \right) \frac{\cplx{v}_k}{\cplx{t}_{\ell}} - \cplx{y}_\ell \cplx{v}_m\right]^*\nonumber\\
&= p_{f\ell} + \jc q_{f\ell}\;\forall \ell = (k,m)\in\Branch,\label{eq1:sf}\\
\cplx{v}_m &\left[- \cplx{y}_\ell \frac{\cplx{v}_k}{\cplx{t}_{\ell}} + \left(\jc \frac{b'_{\ell}}{2} + \cplx{y}_\ell \right) \cplx{v}_m\right]^*\nonumber\\
&= p_{t\ell} + \jc q_{t\ell}\;\forall \ell = (k,m)\in\Branch,\label{eq1:st}
\end{align}
\item Generator power capacities:
\begin{equation}\label{eq1:genlim}
\pl_{Gg}\le p_{Gg} \le \pu_{Gg},\, \ql_{Gg}\le q_{Gg} \le \qu_{Gg}\;\forall g\in\Gen,
\end{equation}
\item Line thermal limits:
\begin{equation}\label{eq1:linelim}
|p_{f\ell} + \jc q_{f\ell}| \le \overline{s}_{\ell},\, |p_{t\ell} + \jc q_{t\ell}| \le \overline{s}_{\ell}\;\forall \ell \in\Branch,
\end{equation}
\item Voltage magnitude limits:
\begin{equation}\label{eq1:magbus}
\vl_k\le \left|\cplx{v}_k\right| \le \vu_k\;\forall k\in\Node,
\end{equation}
\item Reference bus constraint:
\begin{equation}\label{eq1:slack}
\angle \cplx{v}_1 = 0.
\end{equation}
\end{itemize}
\end{subequations}
The objective function \eqref{eq1:objfun} is the cost of conventional generation as commonly used in the literature. Constraints \eqref{eq1:kclp}--\eqref{eq1:st} are derived from Kirchhoff's laws and represent power flows in the network. Constraint~\eqref{eq1:slack} specifies bus $k=1$ as the reference bus. We assume that $\vl_k > 0$ for all $k\in \Node$ in~\eqref{eq1:magbus}, and that  generation cost $c_{g2} p_{Gg}^2 + c_{g1} p_{Gg} + c_{g0}$ is a convex function for all $g\in\Gen$.

Due to the nonconvex constraints~\eqref{eq1:sf}--\eqref{eq1:st},~\eqref{eq1:acopf} is highly nonconvex and NP-hard~\cite{verma2010, lehmann2016}. Applying local methods to this problem provides no guarantee on the optimality of any solution found. Moreover, it is intractable to solve to global optimality for large-scale instances.

\section{ACOPF: Conic relaxations}\label{sec1:conic}
\subsection{Semidefinite relaxation}
With $\cplx{V} = \cv{v}\cv{v}^H$, the ACOPF problem~\eqref{eq1:acopf} can be reformulated as follows
\begin{subequations}\label{eq1:acopfV}
\begin{align}
\min\; &\eqref{eq1:objfun} \nonumber\\
\subj\; &\eqref{eq1:genlim},\eqref{eq1:linelim},\eqref{eq1:slack}, \nonumber\\
&\sum_{g\in \Gen_k} p_{Gg} - p_{Dk} - g_k'\cplx{V}_{kk} = \nonumber\\
&\sum_{\ell=(k,m)\in\Branch} p_{f\ell} + \sum_{\ell=(m,k)\in\Branch} p_{t\ell}\;\forall k\in\Node, \label{eq1:kclpV}\\
&\sum_{g\in \Gen_k} q_{Gg} - q_{Dk} + b_k'\cplx{V}_{kk} = \nonumber\\ &\sum_{\ell=(k,m)\in\Branch} q_{f\ell} + \sum_{\ell=(m,k)\in\Branch} q_{t\ell}\;\forall k\in\Node,\label{eq1:kclqV}\\
&\frac{1}{|\cplx{t}_{\ell}|^2} \left(-\jc \frac{b'_{\ell}}{2} + \cplx{y}_\ell^* \right) \cplx{V}_{kk} - \frac{\cplx{y}_\ell^*}{\cplx{t}_{\ell}} \cplx{V}_{km} \nonumber\\
&= p_{f\ell} + \jc q_{f\ell}\;\forall \ell = (k,m)\in\Branch,\label{eq1:sfV}\\
&-\frac{\cplx{y}_\ell^*}{\cplx{t}_{\ell}^*} \cplx{V}_{mk} + \left(-\jc \frac{b'_{\ell}}{2} + \cplx{y}_\ell^* \right) \cplx{V}_{mm} \nonumber\\
&= p_{t\ell} + \jc q_{t\ell}\;\forall \ell = (k,m)\in\Branch,\label{eq1:stV}\\
&\vl_k^2\le \cplx{V}_{kk} \le \vu_k^2\;\forall k\in\Node,\label{eq1:Vkk}\\
&\cplx{V} = \cv{v}\cv{v}^H.\label{eq1:Vv}
\end{align}
\end{subequations}
The nonconvexity of~\eqref{eq1:acopfV} is captured by the constraint~\eqref{eq1:Vv}. We can show that $\cplx{V} = \cv{v} \cv{v}^H$ if and only if $\cplx{V}\succeq 0$ and $\rank (\cplx{V}) =~ 1$. The \emph{semidefinite relaxation} (SDR) in Model~\ref{model1:sdr} is obtained by dropping the rank constraint. It was first introduced in~\cite{bai} and later, a dual relaxation was developed in~\cite{lavlow3}.
\begin{algorithm}
\caption{Semidefinite relaxation (SDR)}
\label{model1:sdr}
\begin{algorithmic}
\STATE Variables:
\begin{subequations}\label{eq1:var}
\begin{align}
\rv{p}_G, \rv{q}_G &\in \R^{|\Gen|},\\
\rv{p}_f, \rv{q}_f, \rv{p}_t, \rv{q}_t &\in\R^{|\Branch|},\\
\cplx{V} & \in\Herm^{|\Node|}.
\end{align}
\end{subequations}
\STATE Minimize:~\eqref{eq1:objfun}
\STATE Subject to:~\eqref{eq1:genlim},~\eqref{eq1:linelim},~\eqref{eq1:kclpV}--\eqref{eq1:Vkk}, $\cplx{V}\succeq 0$.
\end{algorithmic}
\end{algorithm}

If the optimal solution $\hat{\cplx{V}}$ of SDR is a rank-one matrix, then there exists a complex vector $\hat{\cv{v}}$, global optimal solution of~\eqref{eq1:acopf}. In the literature, there are numerous examples where SDR is exact. However, its exactness is only guaranteed for a few classes of problems under some assumptions~\cite{low2}.

On the other hand, solving SDR for large-scale power systems (more than a thousand of buses) is computationally very expensive. In order to reduce data storage and increase computational speed,~\cite{jabr4} proposes to exploit in SDR the sparsity of the OPF problem. This methodology, as we explain in Section~\ref{subsec1:chord}, suggests to replace the positive semidefinite matrix~$\cplx{V}$ by less-sized positive semidefinite submatrices defined on a chordal extension of the power network~\cite{low1, jabr4, kocuk2016}.

\subsection{Chordal relaxation} \label{subsec1:chord}
Let us interpret the network $\Reseau = (\Node, \Branch)$ as a connected, simple and undirected graph $\G = (\Node,\E)$ where $\Node = \{1,\ldots,n\}$ represents the set of vertices and $\E =\{\{k,m\}: (k,m) \text{ or } (m,k) \in \Branch\}$, the set of edges. 
The power flow equations~\eqref{eq1:kclp}--\eqref{eq1:st} in~\eqref{eq1:acopf} depend only on $\cplx{V}_{kk} := |\cplx{v}_k|^2$, $k\in\Node$, and $\cplx{V}_{km} := \cplx{v}_k \cplx{v}_m^*$, $\{k,m\}\in\E$. In other words, except for the constraint $\cplx{V} \succeq 0$, SDR depends only on a \emph{partial matrix} $\cplx{V}_\G$. A partial matrix means a matrix in which only some of the entries are specified~\cite{grone, low1}.

A subset $\K \subseteq \Node$ is a \emph{clique} if every two distinct vertices in $\K$ are adjacent in $\G$. A clique $\K$ is \emph{maximal} in $\G$ if it is not a subset of a larger clique~$\K'$. A \emph{cycle} is a sequence $k_1 - k_2 - \ldots - k_\gamma - k_1$ of 
$\gamma$ distinct vertices such that $\{k_1,k_2\}, \{k_2,k_3\},\ldots,\{k_{\gamma-1},k_\gamma\},\{k_\gamma,k_1\}\in\E$, where $\gamma \ge 3$ is the length of the cycle. A \emph{chord} of a cycle $k_1 - k_2 - \ldots - k_\gamma - k_1$ is an edge $\{k_i,k_j\} \in \E$ such that $1 \le i < j \le \gamma$ and $2\le j - i \le \gamma-2$. 

$\G$ is \emph{chordal} if every cycle of 4 and more vertices has a chord. A \emph{chordal extension} of $\G$ is a chordal graph $\G' = (\Node,\E')$ that contains $\G$, i.e. $\E \subseteq \E'$. It was proved in~\cite{grone} that the constraint $\cplx{V}\succeq 0$ in SDR is equivalent to $\cplx{V}_\K \succeq 0$ for every maximal clique $\K$ of a chordal extension $\G'$ of $\G$. $\cplx{V}_\K$ is the submatrix of $\cplx{V}$ in which the set of row indices that remain and the set of column indices that remain are both $\K$. Thus, the \emph{chordal relaxation} (CHR) is given in Model~\ref{model1:chr}.

\begin{algorithm}
\caption{Chordal relaxation (CHR)}
\label{model1:chr}
\begin{algorithmic}
\STATE Initialization: $\G = (\Node,\E)$, graph corresponding to $\Reseau = (\Node, \Branch)$. Consider $A = L_\G + I_{|\Node|} \succ 0$, where $L_\G$ is the Laplacian matrix of $\G$ and and $I_{|\Node|}$ is the identity matrix of size $|\Node|$.
\STATE Chordal extension:
\begin{enumerate}
\item Order nodes with heuristic algorithm ``\emph{approximate minimum degree}'' provided by MATLAB-function \texttt{amd}.
\item Compute Cholesky decomposition $LL^T$ of $A$. The sparsity pattern of $L$ defines a chordal extension $\G'$ of $\G$.
\item Identify $\{\mathcal{K}_1,\mathcal{K}_2,\ldots,\mathcal{K}_\kappa\}$, family of maximal cliques of $\G'$.
\end{enumerate}
\STATE Variables:~\eqref{eq1:var}.
\STATE Minimize:~\eqref{eq1:objfun}
\STATE Subject to:~\eqref{eq1:genlim},~\eqref{eq1:linelim},~\eqref{eq1:kclpV}--\eqref{eq1:Vkk}, $\cplx{V}_{\mathcal{K}_i} \succeq 0$, $i = 1,2,\ldots,\kappa$.
\end{algorithmic}
\end{algorithm}

The optimal value $\hat{\upsilon}_{CHR}$ of CHR is not affected by the choice of the chordal extension $\G'$. However, the optimal choice that minimizes the complexity of CHR is NP-hard to compute. Given a positive definite real matrix $A$ of size $n$ such that $A_{km} = 0$ if $\{k,m\}\notin\E$, let $A = LL^T$ be its Cholesky decomposition, where $L$ is a lower triangular matrix. A chordal extension $\G' = (\Node,\E')$ of $\G = (\Node,\E)$ is defined by $\E' = \{\{k,m\}: L_{km} + L_{mk} \ne 0, k\ne m\}$. The fill-in in the Cholesky decomposition depends on the ordering of the nodes $k\in\Node$. The problem of finding the ordering that corresponds to the minimum fill-in is known to be NP-complete. See~\cite{low1,fukuda2001,nakata2003} for more details. Besides,~\cite{jabr4, molzahn2013, andersen2014, madani2014} have developed effective techniques to solve the chordal relaxation of the ACOPF problem and we observe a significant speed-up factor computationally for large-scale power systems compared to the standard SDP relaxation.

\subsection{Second-order cone relaxation}
If we relax the constraint $\cplx{V} \succeq 0$ in SDR by $|\Branch|$ constraints of the form
\begin{equation}\label{eq1:socr}
\cplx{V}_{\{k,m\}} :=
\begin{bmatrix}
\cplx{V}_{kk} & \cplx{V}_{km}\\
\cplx{V}_{km}^* & \cplx{V}_{mm}
\end{bmatrix}
\succeq 0\; \forall (k,m) \in \Branch,\\
\end{equation}
we obtain the standard \emph{second-order cone relaxation} (SOCR) in Model~\ref{model1:socr}. In fact,~\eqref{eq1:socr} represents a rotated second-order cone constraint in the $(\re(\cplx{V}_{km}), \im(\cplx{V}_{km}), \cplx{V}_{kk}, \cplx{V}_{mm})$-space for each branch $(k,m) \in \Branch$.
\begin{algorithm}
\caption{Second-order cone relaxation (SOCR)}
\label{model1:socr}
\begin{algorithmic}
\STATE Variables:~\eqref{eq1:var}.
\STATE Minimize:~\eqref{eq1:objfun}
\STATE Subject to:~\eqref{eq1:genlim},~\eqref{eq1:linelim},~\eqref{eq1:kclpV}--\eqref{eq1:Vkk},~\eqref{eq1:socr}.
\end{algorithmic}
\end{algorithm}

\begin{proposition}[\cite{low1}]\label{prop1:conic}
Let $\hat{\upsilon}$, $\hat{\upsilon}_{SDR}$, $\hat{\upsilon}_{CHR}$, $\hat{\upsilon}_{SOCR}$ be the optimal values of ACOPF Problem~\eqref{eq1:acopf}, SDR, CHR and SOCR. Then $\hat{\upsilon}_{SOCR} \le \hat{\upsilon}_{CHR} = \hat{\upsilon}_{SDR} \le \hat{\upsilon}$. Moreover, for radial networks, $\hat{\upsilon}_{SOCR} = \hat{\upsilon}_{CHR} = \hat{\upsilon}_{SDR} \le \hat{\upsilon}$.
\end{proposition}

SOCR is of significant interest because it is computationally more efficient than SDR, and is thus more amenable for large-scale instances. It was first proposed in~\cite{jabr1} for radial networks, and was extended in~\cite{jabr2} to meshed networks by including a trigonometric functional constraint for the voltage angle spread on each line in the network. Later,~\cite{kocuk2016} proposed three strong SOCP relaxations and showed their computational advantages over SDR.

\section{New conic relaxation}\label{sec1:new}
For two real variables $x$, $y$ such that $\lb{x} \leq x \leq \ub{x}$, $\lb{y} \leq y \leq \ub{y}$  where $\lb{x},\ub{x}, \lb{y}, \ub{y} \in \R$ and $\lb{x}< \ub{x}$, $\lb{y}< \ub{y}$, if $z = xy$ then
\begin{subequations}\label{eq1:rlt}
\begin{align}
z &\leq x\lb{y}+\ub{x}y-\ub{x}\lb{y},\\
z &\leq x\ub{y}+\lb{x}y-\lb{x}\ub{y},\\
z &\geq x\lb{y}+\lb{x}y-\lb{x}\lb{y},\\
z &\geq x\ub{y}+\ub{x}y-\ub{x}\ub{y}.
\end{align}
\end{subequations}
Inequalities~\eqref{eq1:rlt} are called reformulation-linearization technique (RLT) inequalities. They describe the convex hull of $\{(x,y,z)\in \R^3 \colon \lb{x} \leq x \leq \ub{x}, \lb{y} \leq y \leq \ub{y}, z=xy\}$~\cite{mccormick}. For a general nonconvex QCQP with bounded real variables, it has been shown in~\cite{anst} that the use of SDP and RLT constraints together can produce better optimal bounds than either technique used alone. Earlier, it has been proven in~\cite{anstbur} that the convex hull of $\{(\rv{x},X) \in \R^2 \times \Herm^2 \colon X= \rv{x}\rv{x}^T, \lb{\rv{x}} \le~ \rv{x} \le \ub{\rv{x}}\}$ is given by the SDP constraint $X \succeq \rv{x}\rv{x}^T$ together with the RLT inequalities on $X_{11}$, $X_{12}$, $X_{22}$. Hence, one might be tempted to transform ACOPF Problem~\eqref{eq1:acopf} with complex variables $\cv{v}$ into a problem with real variables $\rv{v}^r := \re (\cv{v})$, $\rv{v}^i := \im (\cv{v})$ and consider a relaxation based on SDP and RLT. Such a relaxation would not be as effective as might be expected due to nonrectangular bounds on $\cv{v}$~\cite{chen2017}.

On the other hand, it has been shown in~\cite{josz2016} that relaxing nonconvex constraints of the ACOPF problem before converting from complex to real variables is more advantageous than doing the operations in opposite order. Thus, assuming $|\angle \cplx{v}_k - \angle \cplx{v}_m| \le \pi/2$ for every branch $(k,m)\in\Branch$ in the network, equivalent valid inequalities have been proposed in~\cite{chen2017,coffrin2017} to strengthen the SDP relaxation.

\subsection{Tight-and-cheap relaxation}
For all $(k,m) \in\Branch$, we have $\cplx{V}_{km} = \cplx{v}_k \cplx{v}_m^*$ from \eqref{eq1:Vv}, therefore $|\cplx{V}_{km}| = |\cplx{v}_k| |\cplx{v}_m|$. 
Considering $x = |\cplx{v}_k|$, $y = |\cplx{v}_m|$, $z = |\cplx{V}_{km}|$, and applying \eqref{eq1:rlt}, we obtain
\begin{subequations}\label{eq1:rltVkm}
\begin{gather}
|\cplx{V}_{km}| \leq |\cplx{v}_k| \underline{v}_m + \overline{v}_k |\cplx{v}_m| -\overline{v}_k \underline{v}_m, \label{rltv1}\\
|\cplx{V}_{km}| \leq |\cplx{v}_k| \overline{v}_m + \underline{v}_k |\cplx{v}_m| -\underline{v}_k\overline{v}_m,\\
|\cplx{V}_{km}| \geq |\cplx{v}_k| \underline{v}_m + \underline{v}_k |\cplx{v}_m| -\underline{v}_k \underline{v}_m,\\
|\cplx{V}_{km}| \geq |\cplx{v}_k| \overline{v}_m +\overline{v}_k |\cplx{v}_m| - \overline{v}_k \overline{v}_m, 
\end{gather}
since $\overline{v}_k \le |\cplx{v}_k|\le \overline{v}_k$ for all $k \in \Node$. Moreover for all $k \in\Node$, $\cplx{V}_{kk} := |\cplx{v}_k|^2$, and thus we also have
\begin{equation}\label{eq1:rltVkk}
\cplx{V}_{kk} \le (\vl_k + \vu_k) \left|\cplx{v}_k\right| - \vl_k \vu_k
\end{equation}
\end{subequations}

All RLT inequalities~\eqref{eq1:rltVkm} are nonconvex, except the constraint~\eqref{eq1:rltVkk} corresponding to the reference bus~$k=1$:
\begin{subequations}\label{eq1:rltV11}
\begin{align}
\re (\cplx{v}_1) &\ge \frac{\cplx{V}_{11} + \vl_1 \vu_1}{\vl_1 + \vu_1},\\
\im (\cplx{v}_1) &= 0.
\end{align}
\end{subequations}
Therefore, we define a new formulation of the SDP relaxation in Model~\ref{model1:nsdr}. We denote nSDR. To the best of our knowledge, it is the first time that nSDR with $\cplx{V} \succeq \cv{v}\cv{v}^H$ is proposed for the ACOPF problem.

\begin{algorithm}
\caption{New semidefinite relaxation (nSDR)}
\label{model1:nsdr}
\begin{algorithmic}
\STATE Variables:~\eqref{eq1:var}, $\cv{v} \in \C^{|\Node|}$.
\STATE Minimize:~\eqref{eq1:objfun}
\STATE Subject to:~\eqref{eq1:genlim},~\eqref{eq1:linelim},~\eqref{eq1:kclpV}--\eqref{eq1:Vkk},~\eqref{eq1:rltV11}, $\cplx{V}\succeq \cv{v}\cv{v}^H$.
\end{algorithmic}
\end{algorithm}

\begin{lemma}\label{lemma1:2y}
Let $y \in \R$ such that $\ell \leq y \leq u$, where $0 \le \ell < u < +\infty$. If $x = \sqrt{y}$, then $x \ge \frac{y +  \sqrt{\ell} \sqrt{u}}{\sqrt{\ell} + \sqrt{u}}$.
\end{lemma}
\begin{IEEEproof}
Let $g(y) = \sqrt{y}$ a concave function on its domain. For all $0 \le \ell < u < +\infty$, $\alpha \in [0,1]$,
\[
\begin{aligned}
g((1-\alpha)\ell + \alpha u) &\ge (1-\alpha)g(\ell) + \alpha g(u)\\
&= (1-\alpha)\sqrt{\ell} + \alpha \sqrt{u}.
\end{aligned}
\]
In particular, when $\alpha = \frac{y-\ell}{u-\ell}$, $\ell \leq y \leq u$, we have
\[
\sqrt{y} = x \ge \frac{y +  \sqrt{\ell} \sqrt{u}}{\sqrt{\ell} + \sqrt{u}}.
\]
\end{IEEEproof}
\begin{lemma}\label{lemma1:Vsdp}
Let $\cplx{A} \in \Herm^m$, $\cplx{B} \in \C^{m \times n}$. If $\cplx{A} \succeq 0$, then $\cplx{B}^H \cplx{A} \cplx{B} \succeq~ 0$.
\end{lemma}
\begin{IEEEproof}
Let $\cv{x} \in \C^n$ and $\cv{y} = \cplx{B} \cv{x} \in \C^m$. Therefore, $\cv{x}^H \cplx{B}^H \cplx{A} \cplx{B} \cv{x} = \cv{y}^H \cplx{A} \cv{y} \ge 0$.
\end{IEEEproof}
\begin{proposition}
nSDR is equivalent to SDR.
\end{proposition}
\begin{IEEEproof}
Every feasible solution $\cplx{V}$ of nSDR is also feasible for SDR because $\cplx{V} \succeq~ \cv{v}\cv{v}^H \succeq 0$. 
It remains to prove that for every feasible solution $\cplx{V}$ of SDR, there exists $\cv{v} \in \C^n$ such that $(\cv{v},\cplx{V})$ is feasible for nSDR.\\
Given $\cplx{V}$ feasible solution of SDR, let $\cv{v} = \frac{1}{\sqrt{\cplx{V}_{11}}} \cplx{V} \rv{e}_1$ where $\rv{e}_1$ is the $n$-dimensional vector with $1$ in the first entry and $0$ elsewhere. For all $k \in  \Node$, $ \cplx{v}_k = \frac{1}{\sqrt{\cplx{V}_{11}}} \cplx{V}_{k1}$. In particular, $\cplx{v}_1 = \sqrt{\cplx{V}_{11}} \in \R$ and from Lemma~\ref{lemma1:2y},
\[
\re(\cplx{v}_1) = \sqrt{\cplx{V}_{11}} \ge \frac{\cplx{V}_{11} +  \vl_1 \vu_1}{\vl_1 + \vu_1}.
\]
Now, let $B = \begin{bmatrix}
\frac{\rv{e}_1}{\sqrt{\cplx{V}_{11}}} & I_n
\end{bmatrix} \in \R^{n\times (n+1)}$, where $I_n$ is the identity matrix of size $n$. From Lemma~\ref{lemma1:Vsdp},
\[
B^T \cplx{V} B = \begin{bmatrix}
1 & \cv{v}^H\\
\cv{v} & \cplx{V}
\end{bmatrix} \succeq 0 \Leftrightarrow \cplx{V} \succeq \cv{v}\cv{v}^H.
\]
\end{IEEEproof}

Recall that SOCR is obtained from SDR by replacing the constraint $\cplx{V} \succeq 0$ in SDR by $|\Branch|$ smaller positive semidefiniteness constraints, each one corresponding to a branch of the network. Now we replace the constraint $\cplx{V} \succeq~\cv{v}\cv{v}^H$ in Model~\ref{model1:nsdr} by $|\Branch|$ constraints of the form
\begin{equation}\label{eq1:tcr}
\begin{bmatrix}
1 & \cplx{v}_{k}^* & \cplx{v}_{m}^*\\
\cplx{v}_{k} & \cplx{V}_{kk} & \cplx{V}_{km}\\
\cplx{v}_{m} & \cplx{V}_{km}^* & \cplx{V}_{mm}
\end{bmatrix} \succeq 0\; \forall (k,m)\in\Branch.
\end{equation}
to obtain the relaxation given in Model~\ref{model1:tcr}. We will refer to this relaxation as ``\emph{tight-and-cheap relaxation}'' (TCR). 
Clearly TCR dominates SOCR and is dominated by SDR.
\begin{algorithm}
\caption{Tight-and-cheap relaxation (TCR)}
\label{model1:tcr}
\begin{algorithmic}
\STATE Variables:~\eqref{eq1:var}, $\cv{v}~\in~\C^{|\Node|}$.
\STATE Minimize:~\eqref{eq1:objfun}
\STATE Subject to:~\eqref{eq1:genlim},~\eqref{eq1:linelim},~\eqref{eq1:kclpV}--\eqref{eq1:Vkk},~\eqref{eq1:rltV11},~\eqref{eq1:tcr}.
\end{algorithmic}
\end{algorithm}

\subsection{Strengthening}
Considering $k=1$ as the reference bus, if the constraint $\cplx{V} \succeq 0$ in SDR holds, then
\begin{equation}\label{eq1:stcr}
\cplx{V}_{\{1,k,m\}} := \begin{bmatrix}
\cplx{V}_{11} & \cplx{V}_{1k} & \cplx{V}_{1m}\\
\cplx{V}_{1k}^* & \cplx{V}_{kk} & \cplx{V}_{km}\\
\cplx{V}_{1m}^* & \cplx{V}_{km}^* & \cplx{V}_{mm}
\end{bmatrix} \succeq 0\; \forall (k,m)\in\Branch.
\end{equation}
We define another relaxation given in Model~\ref{model1:stcr}. We call this relaxation ``\emph{strong tight-and-cheap relaxation}'' (STCR). 
Like TCR, STCR dominates SOCR and is dominated by SDR.
\begin{algorithm}
\caption{Strong tight-and-cheap relaxation (STCR)}
\label{model1:stcr}
\begin{algorithmic}
\begin{subequations}
\STATE Variables:~\eqref{eq1:var}.
\STATE Minimize:~\eqref{eq1:objfun}
\STATE Subject to:~\eqref{eq1:genlim},~\eqref{eq1:linelim},~\eqref{eq1:kclpV}--\eqref{eq1:Vkk},~\eqref{eq1:stcr}.
\end{subequations}
\end{algorithmic}
\end{algorithm}
\begin{proposition}\label{prop1:stcr}
STCR is stronger than TCR.
\end{proposition}
\begin{IEEEproof}
We show that for every $\cplx{V}$ feasible solution of STCR, there exists $\cv{v} \in \C^n$ such that $(\cv{v},\cplx{V})$ is feasible for TCR. For all $k \in \Node$, let $ \cplx{v}_k = \frac{1}{\sqrt{\cplx{V}_{11}}} \cplx{V}_{k1}$. In particular, $\cplx{v}_1 = \sqrt{\cplx{V}_{11}} \in \R$ and from Lemma~\ref{lemma1:2y},
\[
\re(\cplx{v}_1) = \sqrt{\cplx{V}_{11}} \ge \frac{\cplx{V}_{11} +  \vl_1 \vu_1}{\vl_1 + \vu_1}.
\]
Now, for all $(k,m) \in \Branch$, if~\eqref{eq1:stcr} holds, then by Lemma~\ref{lemma1:Vsdp},
\begin{gather*}
\begin{bmatrix}
\frac{1}{\sqrt{\cplx{V}_{11}}} & \rv{0}^T\\
\rv{0} & I_2\\
\end{bmatrix}
\cplx{V}_{\{1,k,m\}}
\begin{bmatrix}
\frac{1}{\sqrt{\cplx{V}_{11}}} & \rv{0}^T\\
\rv{0} & I_2\\
\end{bmatrix}\\
=
\begin{bmatrix}
1 & \cplx{v}_{k}^* & \cplx{v}_{m}^*\\
\cplx{v}_{k} & \cplx{V}_{kk} & \cplx{V}_{km}\\
\cplx{v}_{m} & \cplx{V}_{km}^* & \cplx{V}_{mm}
\end{bmatrix}
\succeq 0.
\end{gather*}
\end{IEEEproof}
\begin{proposition}
Let $\G = (\Node,\E)$ a graph corresponding to a power network $\Reseau = (\Node,\Branch)$. Suppose bus $k=1$ is the reference bus. If the induced subgraph $\G - \{1\}$ has no cycle, then STCR is equivalent to SDR.
\end{proposition}
\begin{IEEEproof}
Consider $\G' = (\Node,\E')$ where $\E' = \E \cup \{\{1,m\}\colon \{1,m\}\notin \E\}$. Since $\G - \{1\}$ has no cycle, every cycle in $\G'$ contains vertex $k=1$. Also, since for all $m\in\Node\setminus\{1\}$, $\{1,m\}\in\E'$, every cycle $1-m_1-m_2-m_3-1$ of 4 vertices has a chord $\{1,m_2\}$. Thus, $\G'$ is a chordal extension of~$\G$.\\
On the other hand, since $\G - \{1\}$ has no cycle, every maximal clique $\K$ of $\G'$ contains vertex $k=1$ and has at most 3 vertices, i.e. $\K = \{1,k,m\}$ for all $\{k,m\}\in\E$. Then the constraint $\cplx{V}\succeq 0$ in SDR is equivalent to $\cplx{V}_{\{1,k,m\}} \succeq 0$ for all $\{k,m\} \in \E$.
\end{IEEEproof}

\section{Computational results}\label{sec1:results}
In this section, we evaluate the accuracy and the computational efficiency of TCR and STCR as compared to SOCR, CHR and SDR.

We tested the models~\ref{model1:sdr},~\ref{model1:chr},~\ref{model1:socr},~\ref{model1:tcr} and~\ref{model1:stcr} on standard test cases available from MATPOWER~\cite{matpower, birch, josz2016ac}. {\em It is important to note that, unlike what was
done in~\cite{lavlow3, molzahn2013, andersen2014, madani2014}, we did not make any modification or simplification to the data.} 

We solved all the relaxations in MATLAB using CVX~2.1~\cite{cvx2, gb08} with the solver MOSEK~8.0.0.60 and \texttt{default precision} (tolerance $\epsilon = 1.49 \times 10^{-8}$). All the computations were carried out on an \texttt{Intel Core i7-6700 CPU @ 3.40 GHz} computing platform. When solving SOCR for instances with at least 1000 buses, MOSEK ended its computation with message \texttt{Mosek error: MSK\_RES\_TRM\_STALL()}. For these test cases, we replace
constraints~\eqref{eq1:socr} by equivalent ones
\[
\begin{bmatrix}
\cplx{V}_{kk} + \cplx{V}_{mm} & 0 & 2\cplx{V}_{km}\\
0 & \cplx{V}_{kk} + \cplx{V}_{mm} & \cplx{V}_{kk} - \cplx{V}_{mm}\\
2\cplx{V}_{km}^*& \cplx{V}_{kk} - \cplx{V}_{mm} & \cplx{V}_{kk} + \cplx{V}_{mm}
\end{bmatrix} \succeq 0
\]
for all $(k,m)\in\Branch$.

We considered two objective functions: the generation cost~[\$/h]~\eqref{eq1:objfun} and the active loss [MW] where $c_{g2} = 0$, $c_{g1} = 1$ and $c_{g0} = 0$ for all $g\in\Gen$ in~\eqref{eq1:objfun}. Both objective functions of test cases from~\cite{josz2016ac} are the same. We denote $\lb{\upsilon}$ the best lower bound which is the maximum value among $\hat{\upsilon}_{SOCR}$, $\hat{\upsilon}_{TCR}$, $\hat{\upsilon}_{STCR}$, $\hat{\upsilon}_{CHR}$, $\hat{\upsilon}_{SDR}$, respective optimal values of SOCR, TCR, STCR, CHR and SDR. The optimality gap is measured as $100(1-\hat{\upsilon}_R/\ub{\upsilon})$, where $\ub{\upsilon}$ is the upper bound provided by the MATPOWER-solver ``\texttt{MIPS}'' and $\hat{\upsilon}_R$ is the relaxation optimal value. For some test cases: \texttt{1888rte}, \texttt{1951rte}, \texttt{2848rte}, \texttt{2868rte}, \texttt{6468rte}, \texttt{6470rte}, \texttt{6495rte} and \texttt{6515rte}, \texttt{MIPS} failed to find a local optimal solution, so we considered the upper bounds reported in~\cite{josz2016ac}. 

Table~\ref{table1:cost} and Table~\ref{table1:loss} summarize the optimality gaps of the five relaxations for cost minimization and loss minimization, respectively. The results support the following key points: 
\begin{enumerate}
\item CHR is equivalent to SDR as predicted by Proposition~\ref{prop1:conic}.
\item TCR and STCR are stronger than SOCR. When compared to SOCR, TCR reduces the optimality gap from 0.17\% to 0.06\% on average for large-scale instances in Table~\ref{table1:loss}.
\item Optimality gaps of TCR and STCR are very close to CHR or SDR. We observe significant optimality gaps of TCR and STCR when the optimality gap of CHR or SDR is not close to zero, e.g. \texttt{case5} and \texttt{case\_ACTIV\_SG\_500} instances in Table~\ref{table1:cost}.
\item STCR is stronger than TCR as predicted by Proposition~\ref{prop1:stcr}. For example, STCR reduces substantially the optimality gap of TCR from 12.75\% to 5.22\% on \texttt{case5} instance in Table~\ref{table1:cost}.
\end{enumerate}

The computation times reported by MOSEK are shown in Table~\ref{table1:cost} and Table~\ref{table1:loss}. The time CVX took to pre-compile a model is not included. For CHR, the computation time does not take into account the time of building the chordal extension of an instance's graph. We did not solve SDR for the extra large-scale instances (those with at least 6 000 buses) because of the high computational cost. In Tables~\ref{table1:cost}~and~\ref{table1:loss}, we note:
\begin{enumerate}
\item Among all relaxations, SOCR is the fastest and SDR is the slowest.
\item CHR is on average around 30 times faster than SDR for large-scale instances.
\item TCR is on average around 30 times faster than CHR for large-scale instances and 55 times for extra large-scale instances.
\item TCR is on average around 3 times faster than STCR for large-scale instances and 7 times for extra large-scale instances.
\end{enumerate}

\section{Conclusion}\label{sec1:conclusion}
We proposed a new formulation of the semidefinite relaxation for the ACOPF problem. This formulation is based on a positive semidefiniteness constraint combined with reformulation-linearization technique (RLT) constraints defined on the reference (slack) bus. We proved that it is equivalent to the standard SDP relaxation. Thereafter, we derived a tight-and-cheap semidefinite relaxation (TCR) stronger than the standard SOCP relaxation. Experiments on unmodified MATPOWER instances show that the proposed relaxation offers an interesting trade-off between the standard SDP and SOCP relaxations for large-scale power systems because it is very close to the SDP relaxation in terms of optimality gap, but computationally it is much faster than the chordal relaxation (which is equivalent to the SDP relaxation).

A strong TCR (STCR) was also proposed. We showed that, under some assumption, it is tantamount to the standard SDP relaxation. Although faster than the chordal relaxation, it is not as fast as TCR.

\begin{table*}[!t]
\footnotesize
\centering
\caption{Cost minimization}\label{table1:cost}
\resizebox{\linewidth}{!}{
\begin{tabular}{@{}lrr|rrrrr|rrrrr@{}}
\toprule
Test case & $\ub{\upsilon}$ [\$/h] & $\lb{\upsilon}$ [\$/h] & \multicolumn{5}{c|}{Optimality gap [\%]} & \multicolumn{5}{c}{Computation time [s]} \\
\cmidrule(r){4-8} \cmidrule(l){9-13} & 	&  & SOCR &  TCR & STCR & CHR & SDR & SOCR &  TCR & STCR & CHR & SDR \\
\midrule
\multicolumn{13}{c}{\it Small-scale instances}\\
\midrule
{\tt LMBD3\_50}	&	5 812.64	&	5 789.91	&	1.32	&	0.74	&	0.39	&	0.39	&	0.39	&	0.08	&	0.06	&	0.06	&	0.06	&	0.06	\\
{\tt LMBD3\_60}	&	5 707.11	&	5 707.11	&	0.05	&	0.00	&	0.00	&	0.00	&	0.00	&	0.07	&	0.08	&	0.06	&	0.06	&	0.05	\\
{\tt case5}	&	17 551.89	&	16 635.78	&	14.54	&	12.75	&	5.22	&	5.22	&	5.22	&	0.06	&	0.07	&	0.07	&	0.07	&	0.07	\\
{\tt case6ww}	&	3 143.97	&	3 143.97	&	0.63	&	0.00	&	0.00	&	0.00	&	0.00	&	0.09	&	0.07	&	0.07	&	0.07	&	0.06	\\
{\tt case9}	&	5 296.69	&	5 296.69	&	0.00	&	0.00	&	0.00	&	0.00	&	0.00	&	0.07	&	0.07	&	0.07	&	0.07	&	0.07	\\
{\tt case14}	&	8 081.53	&	8 081.52	&	0.08	&	0.00	&	0.00	&	0.00	&	0.00	&	0.07	&	0.08	&	0.07	&	0.08	&	0.07	\\
{\tt case24\_ieee\_rts}	&	63 352.21	&	63 352.20	&	0.01	&	0.00	&	0.00	&	0.00	&	0.00	&	0.09	&	0.13	&	0.14	&	0.11	&	0.13	\\
{\tt case30}	&	576.89	&	576.89	&	0.57	&	0.07	&	0.00	&	0.00	&	0.00	&	0.08	&	0.14	&	0.16	&	0.13	&	0.16	\\
{\tt case\_ieee30}	&	8 906.14	&	8 906.14	&	0.04	&	0.00	&	0.00	&	0.00	&	0.00	&	0.08	&	0.11	&	0.11	&	0.08	&	0.10	\\
{\tt case39}	&	41 864.18	&	41 862.03	&	0.02	&	0.01	&	0.01	&	0.01	&	0.01	&	0.08	&	0.19	&	0.18	&	0.11	&	0.25	\\
{\tt case57}	&	41 737.79	&	41 737.78	&	0.06	&	0.01	&	0.00	&	0.00	&	0.00	&	0.08	&	0.17	&	0.19	&	0.16	&	0.25	\\
{\tt case89pegase}	&	5 819.81	&	5 819.65	&	0.17	&	0.04	&	0.00	&	0.00	&	0.00	&	0.16	&	0.68	&	0.83	&	0.90	&	1.08	\\
{\bf Average}	&	&	&	\bf 0.39	&	\bf 0.23	&	\bf 0.09	&	\bf 0.09	&	\bf 0.09	&	\bf 0.10	&	\bf 0.30	&	\bf 0.34	&	\bf 0.34	&	\bf 0.43	\\
\midrule
\multicolumn{13}{c}{\it Medium-scale instances}\\
\midrule
{\tt case118}	&	129 660.70	&	129 654.54	&	0.25	&	0.03	&	0.02	&	0.00	&	0.00	&	0.08	&	0.33	&	0.42	&	0.32	&	1.09	\\
{\tt case\_ACTIV\_SG\_200}	&	27 557.57	&	27 557.55	&	0.00	&	0.00	&	0.00	&	0.00	&	0.00	&	0.21	&	0.56	&	0.65	&	0.72	&	4.66	\\
{\tt case\_illinois200}	&	36 748.39	&	36 748.33	&	0.02	&	0.00	&	0.00	&	0.00	&	0.00	&	0.22	&	0.83	&	0.92	&	0.87	&	6.28	\\
{\tt case300}	&	719 725.11	&	719 710.63	&	0.15	&	0.02	&	0.01	&	0.00	&	0.00	&	0.20	&	0.97	&	1.20	&	0.95	&	10.52	\\
{\tt case\_ACTIV\_SG\_500}	&	72 578.30	&	71 048.04	&	5.38	&	4.39	&	4.20	&	2.11	&	2.11	&	0.75	&	3.21	&	4.04	&	3.14	&	103.35	\\
{\bf Average}	&	&	&	\bf 2.10	&	\bf 1.67	&	\bf 1.60	&	\bf 0.80	&	\bf 0.80	&	\bf 0.40	&	\bf 1.68	&	\bf 2.08	&	\bf 1.68	&	\bf 43.36	\\
\midrule
\multicolumn{13}{c}{\it Large-scale instances}\\
\midrule
{\tt case1354pegase}	&	74 069.35	&	74 061.72	&	0.08	&	0.02	&	0.02	&	0.01	&	0.01	&	5.73	&	6.34	&	13.23	&	9.72	&	1 657.85	\\
{\tt case1888rte}	&	59 805.1\phantom{0}	&	59 601.29	&	0.39	&	0.36	&	0.35	&	0.34	&	0.34	&	8.88	&	10.61	&	26.64	&	17.88	&	4 629.31	\\
{\tt case1951rte}	&	81 737.7\phantom{0}	&	81 725.16	&	0.08	&	0.03	&	0.03	&	0.01	&	0.02	&	11.24	&	12.49	&	31.78	&	23.31	&	5 595.59	\\
{\tt case2383wp}	&	1 868 511.83	&	1 861 214.70	&	1.07	&	0.50	&	0.48	&	0.40	&	0.39	&	11.55	&	11.60	&	33.76	&	216.38	&	9 420.67	\\
{\tt case2736sp}	&	1 307 883.13	&	1 307 695.31	&	0.31	&	0.03	&	0.01	&	0.02	&	0.02	&	10.24	&	12.28	&	36.42	&	303.63	&	10 786.69	\\
{\tt case2737sop}	&	777 629.30	&	777 517.52	&	0.27	&	0.03	&	0.02	&	0.02	&	0.01	&	8.90	&	11.47	&	35.68	&	242.92	&	10 601.37	\\
{\tt case2746wop}	&	1 208 279.81	&	1 208 182.18	&	0.40	&	0.03	&	0.02	&	0.01	&	0.01	&	9.87	&	11.51	&	35.57	&	343.90	&	10 480.83	\\
{\tt case2746wp}	&	1 631 775.10	&	1 631 665.81	&	0.33	&	0.03	&	0.02	&	0.01	&	0.01	&	11.34	&	13.00	&	36.68	&	340.18	&	11 654.59	\\
{\tt case2848rte}	&	53 021.8\phantom{0}	&	53 005.20	&	0.08	&	0.04	&	0.04	&	0.03	&	0.04	&	12.17	&	14.36	&	44.56	&	46.64	&	14 567.98	\\
{\tt case2868rte}	&	79 794.7\phantom{0}	&	79 787.67	&	0.07	&	0.02	&	0.02	&	0.01	&	0.02	&	14.05	&	17.03	&	55.08	&	47.20	&	16 933.60	\\
{\tt case2869pegase}	&	133 999.29	&	133 983.11	&	0.09	&	0.03	&	0.03	&	0.01	&	0.03	&	15.45	&	18.76	&	60.18	&	49.29	&	14 866.00	\\
{\tt case3012wp}	&	2 591 706.57	&	2 587 512.73	&	0.82	&	0.38	&	0.37	&	0.17	&	0.16	&	12.90	&	15.21	&	49.39	&	1 276.67	&	16 051.95	\\
{\tt case3120sp}	&	2 142 703.76	&	2 140 385.00	&	0.56	&	0.13	&	0.12	&	0.11	&	0.11	&	14.98	&	16.55	&	52.77	&	765.39	&	14 636.04	\\
{\tt case3375wp}	&	7 412 030.68	&	7 407 116.46	&	0.30	&	0.14	&	0.14	&	0.08	&	0.07	&	15.23	&	17.64	&	66.37	&	1 216.15	&	18 229.90	\\
{\bf Average}	&	&	&	\bf 0.36	&	\bf 0.12	&	\bf 0.12	&	\bf 0.08	&	\bf 0.08	&	\bf 12.03	&	\bf 14.02	&	\bf 43.71	&	\bf 402.09	&	\bf 12 336.09	\\
\midrule
\multicolumn{13}{c}{\it Extra large-scale instances}\\
\midrule
{\tt case6468rte}	&	86 860.0\phantom{0}	&	86 808.12	&	0.27	&	0.08	&	0.08	&	0.06	&	--	&	35.47	&	40.79	&	314.92	&	1 990.53	&	--	\\
{\tt case6470rte}	&	98 345.5\phantom{0}	&	98 333.90	&	0.18	&	0.06	&	0.03	&	0.01	&	--	&	53.34	&	51.88	&	371.42	&	2 729.29	&	--	\\
{\tt case6495rte}	&	106 283.4\phantom{0}	&	106 072.44	&	0.46	&	0.23	&	0.23	&	0.20	&	--	&	50.50	&	65.15	&	431.65	&	3 592.15	&	--	\\
{\tt case6515rte}	&	109 804.2\phantom{0}	&	109 688.71	&	0.38	&	0.16	&	0.14	&	0.11	&	--	&	48.02	&	59.60	&	427.22	&	3 523.01	&	--	\\
{\bf Average}	&	&	&	\bf 0.32	&	\bf 0.13	&	\bf 0.12	&	\bf 0.09	&	\bf --	&	\bf 46.84	&	\bf 54.38	&	\bf 386.42	&	\bf 2 960.41	&	\bf --	\\
\bottomrule
\end{tabular}
}
\end{table*}

\begin{table*}[!t]
\footnotesize
\centering
\caption{Loss minimization}\label{table1:loss}
\resizebox{\linewidth}{!}{
\begin{tabular}{@{}lrr|rrrrr|rrrrr@{}}
\toprule
Test case & $\ub{\upsilon}$ [MW] & $\lb{\upsilon}$ [MW] & \multicolumn{5}{c|}{Optimality gap [\%]}	& \multicolumn{5}{c}{Computation time [s]} \\
\cmidrule(r){4-8} \cmidrule(l){9-13} & 	&  & SOCR &  TCR & STCR & CHR & SDR & SOCR &  TCR & STCR & CHR & SDR \\
\midrule
\multicolumn{13}{c}{\it Small-scale instances}\\
\midrule
{\tt LMBD3\_50}	&	317.38	&	317.38	&	0.00	&	0.00	&	0.00	&	0.00	&	0.00	&		0.05	&	0.06	&	0.06	&	0.05	&	0.06	\\
{\tt LMBD3\_60}	&	316.75	&	316.75	&	0.01	&	0.00	&	0.00	&	0.00	&	0.00	&	0.05	&	0.05	&	0.07	&	0.05	&	0.06	\\
{\tt case5}	&	1 001.06	&	1 001.06	&	0.00	&	0.00	&	0.00	&	0.00	&	0.00	&	0.08	&	0.10	&	0.09	&	0.09	&	0.09	\\
{\tt case6ww}	&	216.84	&	216.84	&	0.16	&	0.00	&	0.00	&	0.00	&	0.00	&	0.06	&	0.07	&	0.06	&	0.08	&	0.06	\\
{\tt case9}	&	317.32	&	317.32	&	0.00	&	0.00	&	0.00	&	0.00	&	0.00	&	0.07	&	0.07	&	0.07	&	0.07	&	0.07	\\
{\tt case14}	&	259.55	&	259.55	&	0.00	&	0.00	&	0.00	&	0.00	&	0.00	&	0.07	&	0.07	&	0.07	&	0.06	&	0.05	\\
{\tt case24\_ieee\_rts}	&	2 875.75	&	2 875.74	&	0.01	&	0.00	&	0.00	&	0.00	&	0.00	&	0.10	&	0.14	&	0.14	&	0.12	&	0.16	\\
{\tt case30}	&	191.09	&	191.09	&	0.23	&	0.01	&	0.00	&	0.00	&	0.00	&	0.08	&	0.13	&	0.15	&	0.10	&	0.18	\\
{\tt case\_ieee30}	&	284.77	&	284.77	&	0.05	&	0.00	&	0.00	&	0.00	&	0.00	&	0.07	&	0.09	&	0.09	&	0.08	&	0.27	\\
{\tt case39}	&	6 284.15	&	6 283.90	&	0.01	&	0.00	&	0.00	&	0.00	&	0.00	&	0.09	&	0.19	&	0.18	&	0.12	&	0.22	\\
{\tt case57}	&	1 262.10	&	1 262.10	&	0.03	&	0.00	&	0.00	&	0.00	&	0.00	&	0.08	&	0.14	&	0.16	&	0.15	&	0.21	\\
{\tt case89pegase}	&	5 819.81	&	5 819.65	&	0.17	&	0.04	&	0.00	&	0.00	&	0.00	&	0.16	&	0.68	&	0.83	&	0.90	&	1.08	\\
{\bf Average}	&	&	&	\bf 0.09	&	\bf 0.01	&	\bf 0.00	&	\bf 0.00	&	\bf 0.00	&	\bf 0.10	&	\bf 0.29	&	\bf 0.33	&	\bf 0.34	&	\bf 0.44	\\
\midrule
\multicolumn{13}{c}{\it Medium-scale instances}\\
\midrule
{\tt case118}	&	4 251.23	&	4 251.03	&	0.01	&	0.01	&	0.00	&	0.00	&	0.00	&	0.09	&	0.31	&	0.35	&	0.25	&	0.83	\\
{\tt case\_ACTIV\_SG\_200}	&	1 483.92	&	1 483.92	&	0.01	&	0.00	&	0.00	&	0.00	&	0.00	&	0.19	&	0.52	&	0.65	&	0.56	&	4.07	\\
{\tt case\_illinois200}	&	2 246.49	&	2 246.48	&	0.01	&	0.00	&	0.00	&	0.00	&	0.00	&	0.25	&	0.79	&	0.90	&	0.83	&	5.80	\\
{\tt case300}	&	23 737.72	&	23 737.55	&	0.06	&	0.01	&	0.01	&	0.00	&	0.00	&	0.19	&	0.91	&	1.05	&	0.86	&	9.93	\\
{\tt case\_ACTIV\_SG\_500}	&	7 817.46	&	7 817.41	&	0.02	&	0.00	&	0.00	&	0.00	&	0.00	&	0.75	&	4.04	&	5.03	&	3.75	&	127.07	\\
{\bf Average}	&	&	&	\bf 0.02	&	\bf 0.00	&	\bf 0.00	&	\bf 0.00	&	\bf 0.00	&	\bf 0.40	&	\bf 1.97	&	\bf 2.41	&	\bf 1.85	&	\bf 52.04	\\
\midrule
\multicolumn{13}{c}{\it Large-scale instances}\\
\midrule
{\tt case1354pegase}	&	74 069.35	&	74 061.72	&	0.08	&	0.02	&	0.02	&	0.01	&	0.01	&	5.73	&	6.34	&	13.23	&	9.72	&	1 657.85	\\
{\tt case1888rte}	&	59 805.1\phantom{0}	&	59 601.29	&	0.39	&	0.36	&	0.35	&	0.34	&	0.34	&	8.88	&	10.61	&	26.64	&	17.88	&	4 629.31	\\
{\tt case1951rte}	&	81 737.7\phantom{0}	&	81 725.16	&	0.08	&	0.03	&	0.03	&	0.01	&	0.02	&	11.24	&	12.49	&	31.78	&	23.31	&	5 595.59	\\
{\tt case2383wp}	&	24 991.40	&	24 979.30	&	0.21	&	0.07	&	0.07	&	0.05	&	0.05	&	10.83	&	10.88	&	32.25	&	188.90	&	7 758.35	\\
{\tt case2736sp}	&	18 335.95	&	18 334.71	&	0.19	&	0.03	&	0.02	&	0.01	&	0.01	&	9.27	&	11.01	&	31.96	&	290.60	&	10 768.88	\\
{\tt case2737sop}	&	11 397.39	&	11 396.61	&	0.18	&	0.02	&	0.01	&	0.01	&	0.01	&	8.46	&	11.38	&	29.65	&	216.13	&	9 715.82	\\
{\tt case2746wop}	&	19 212.35	&	19 211.42	&	0.21	&	0.03	&	0.02	&	0.01	&	0.00	&	10.15	&	11.26	&	30.92	&	320.35	&	9 093.47	\\
{\tt case2746wp}	&	25 269.45	&	25 268.43	&	0.19	&	0.03	&	0.01	&	0.00	&	0.01	&	10.76	&	12.34	&	33.78	&	336.99	&	11 785.35	\\
{\tt case2848rte}	&	53 021.8\phantom{0}	&	53 005.20	&	0.08	&	0.04	&	0.04	&	0.03	&	0.04	&	12.17	&	14.36	&	44.56	&	46.64	&	14 567.98	\\
{\tt case2868rte}	&	79 794.7\phantom{0}	&	79 787.67	&	0.07	&	0.02	&	0.02	&	0.01	&	0.02	&	14.05	&	17.03	&	55.08	&	47.20	&	16 933.60	\\
{\tt case2869pegase}	&	133 999.29	&	133 983.11	&	0.09	&	0.03	&	0.03	&	0.01	&	0.03	&	15.45	&	18.76	&	60.18	&	49.29	&	14 866.00	\\
{\tt case3012wp}	&	27 645.97	&	27 637.25	&	0.22	&	0.05	&	0.04	&	0.04	&	0.03	&	16.80	&	17.57	&	55.82	&	1 448.81	&	17 689.12	\\
{\tt case3120sp}	&	21 513.52	&	21 495.85	&	0.24	&	0.10	&	0.09	&	0.09	&	0.08	&	15.08	&	17.09	&	58.92	&	830.41	&	15 981.61	\\
{\tt case3375wp}	&	49 004.69	&	48 995.70	&	0.15	&	0.04	&	0.03	&	0.02	&	0.03	&	15.07	&	16.85	&	70.58	&	1 229.45	&	17 306.68	\\
{\bf Average}	&	&	&	\bf 0.17	&	\bf 0.06	&	\bf 0.05	&	\bf 0.04	&	\bf 0.04	&	\bf 12.17	&	\bf 13.97	&	\bf 43.70	&	\bf 416.32	&	\bf 12 230.42	\\
\midrule
\multicolumn{13}{c}{\it Extra large-scale instances}\\
\midrule
{\tt case6468rte}	&	86 860.0\phantom{0}	&	86 808.12	&	0.27	&	0.08	&	0.08	&	0.06	&	--	&	35.47	&	40.79	&	314.92	&	1 990.53	&	--	\\
{\tt case6470rte}	&	98 345.5\phantom{0}	&	98 333.90	&	0.18	&	0.06	&	0.03	&	0.01	&	--	&	53.34	&	51.88	&	371.42	&	2 729.29	&	--	\\
{\tt case6495rte}	&	106 283.4\phantom{0}	&	106 072.44	&	0.46	&	0.23	&	0.23	&	0.20	&	--	&	50.50	&	65.15	&	431.65	&	3 592.15	&	--	\\
{\tt case6515rte}	&	109 804.2\phantom{0}	&	109 688.71	&	0.38	&	0.16	&	0.14	&	0.11	&	--	&	48.02	&	59.60	&	427.22	&	3 523.01	&	--	\\
{\bf Average}	&	&	&	\bf 0.32	&	\bf 0.13	&	\bf 0.12	&	\bf 0.09	&	\bf --	&	\bf 46.84	&	\bf 54.38	&	\bf 386.42	&	\bf 2 960.41	&	\bf --	\\
\bottomrule
\end{tabular}
}
\end{table*}

% if have a single appendix:
%\appendix[Proof of the Zonklar Equations]
% or
%\appendix  % for no appendix heading
% do not use \section anymore after \appendix, only \section*
% is possibly needed

% use appendices with more than one appendix
% then use \section to start each appendix
% you must declare a \section before using any
% \subsection or using \label (\appendices by itself
% starts a section numbered zero.)
%

%\appendices
%\section{Proof of the First Zonklar Equation}
%Appendix one text goes here.
%
%% you can choose not to have a title for an appendix
%% if you want by leaving the argument blank
%\section{}
%Appendix two text goes here.

% use section* for acknowledgment
\section*{Acknowledgment}
We thank St\'ephane Alarie and Laurent Lenoir, both of the Hydro-Qu\'ebec Research Institute (IREQ), for helpful comments on early drafts of this paper. We also thank the anonymous reviewers for their many helpful suggestions that helped us improve this paper.

% Can use something like this to put references on a page
% by themselves when using endfloat and the captionsoff option.
\ifCLASSOPTIONcaptionsoff
  \newpage
\fi

% trigger a \newpage just before the given reference
% number - used to balance the columns on the last page
% adjust value as needed - may need to be readjusted if
% the document is modified later
%\IEEEtriggeratref{8}
% The "triggered" command can be changed if desired:
%\IEEEtriggercmd{\enlargethispage{-5in}}

% references section

% can use a bibliography generated by BibTeX as a .bbl file
% BibTeX documentation can be easily obtained at:
% http://mirror.ctan.org/biblio/bibtex/contrib/doc/
% The IEEEtran BibTeX style support page is at:
% http://www.michaelshell.org/tex/ieeetran/bibtex/
%\bibliographystyle{IEEEtran}
% argument is your BibTeX string definitions and bibliography database(s)
%\bibliography{IEEEabrv,../bib/paper}
%
% <OR> manually copy in the resultant .bbl file
% set second argument of \begin to the number of references
% (used to reserve space for the reference number labels box)
\bibliographystyle{unsrt}
\bibliography{../drmth_biblio}
%\begin{thebibliography}{1}
%
%\bibitem{IEEEhowto:kopka}
%H.~Kopka and P.~W. Daly, \emph{A Guide to \LaTeX}, 3rd~ed.\hskip 1em plus
%  0.5em minus 0.4em\relax Harlow, England: Addison-Wesley, 1999.
%
%\end{thebibliography}

% biography section
% 
% If you have an EPS/PDF photo (graphicx package needed) extra braces are
% needed around the contents of the optional argument to biography to prevent
% the LaTeX parser from getting confused when it sees the complicated
% \includegraphics command within an optional argument. (You could create
% your own custom macro containing the \includegraphics command to make things
% simpler here.)
%\begin{IEEEbiography}[{\includegraphics[width=1in,height=1.25in,clip,keepaspectratio]{mshell}}]{Michael Shell}
% or if you just want to reserve a space for a photo:

\begin{IEEEbiographynophoto}{Christian Bingane}
(S'18) received the B.Eng. degree in electrical engineering in 2014 from Polytechnique Montreal, Montreal, QC, Canada, where he is currently working toward the Ph.D. degree in applied mathematics. He is currently a student member of the GERAD research center.

His research interests include optimization in power systems and conic programming. He is concerned with using linear programming, second-order cone programming or semidefinite programming to provide guaranteed global optimal solution to the optimal power flow problem for a large-scale power system.
\end{IEEEbiographynophoto}

% if you will not have a photo at all:
\begin{IEEEbiographynophoto}{Miguel F. Anjos}
(M'07--SM'18) received the B.Sc. degree, the M.S. degree and the Ph.D. degree from McGill University, Montreal, QC, Canada, Stanford University, Stanford, CA, USA, and the University of Waterloo, Waterloo, ON, Canada in 1992, 1994 and 2001 respectively. 

He is currently a Professor with the Department of Mathematics and Industrial Engineering, Polytechnique Montreal, Montreal, QC, Canada, where he holds the NSERC-Hydro-Quebec-Schneider Electric Industrial Research Chair, and an Inria International Chair. He is a Licensed Professional Engineer in Ontario, Canada. He served for five years as Editor-in-Chief of Optimization and Engineering, and serves on several editorial boards.

His allocades include a Canada Research Chair, the M\'eritas Teaching Award, a Humboldt Research Fellowship, the title of EUROPT Fellow, and the Queen Elizabeth II Diamond Jubilee Medal. He is an elected Fellow of the Canadian Academy of Engineering. 
\end{IEEEbiographynophoto}

% insert where needed to balance the two columns on the last page with
% biographies
%\newpage

\begin{IEEEbiographynophoto}{S\'ebastien Le~Digabel}
received the M.Sc.A. degree and the Ph.D. degree in applied mathematics from Polytechnique Montreal, Montreal, Quebec, Canada in 2002 and 2008 respectively. He was a postdoctoral fellow with the IBM Watson Research Center and the University of Chicago in 2010 and 2011.

He is currently an Associate Professor with the Department of Mathematics and Industrial Engineering, Polytechnique Montreal, Montreal, QC, Canada, and a regular member of the GERAD research center.

His research interests include the analysis and development of algorithms for derivative-free and blackbox optimization, and the design of related software. All of his work on derivative-free optimization is included in the NOMAD software, a free package for blackbox optimization available at {\tt www.gerad.ca/nomad}.
\end{IEEEbiographynophoto}

% You can push biographies down or up by placing
% a \vfill before or after them. The appropriate
% use of \vfill depends on what kind of text is
% on the last page and whether or not the columns
% are being equalized.

%\vfill

% Can be used to pull up biographies so that the bottom of the last one
% is flush with the other column.
%\enlargethispage{-5in}

% that's all folks
\end{document}